\documentclass[11pt]{article}
\usepackage{times,amssymb,amsfonts}
\def\qed{\ifhmode\unskip\nobreak\fi\ifmmode\ifinner\else\hskip.5em\fi\fi
 \hbox{\hskip.5em$\square$\hskip.1em}}
\newenvironment{genericem}[1]{\smallskip{\sc #1.}\em}{\/\rm \smallskip}
\newenvironment{proof}{\smallskip{\sc Proof.}}{\qed\smallskip}

\begin{document}
\title{Stable Configurations of repelling Points on compact Manifolds}
\author{Burton Randol}
\date{}
\maketitle

\abstract{\sf \scriptsize\begin{sloppypar}This is an expanded version of \cite{randolarxiv_1}. Using techniques from \cite{randolchapter}, in which a differential-geometrically intrinsic treatment of counterparts of classical electrostatics was introduced, it is shown that on some compact manifolds, certain stable configurations of points which mutually repel along all interconnecting geodesics become equidistributed as the number of points increases.\end{sloppypar}}

\bigskip

\bigskip

In this note, we describe a feature of stable configurations of repelling points on certain manifolds in the sense of \cite{randolchapter}\ (pp.\ 282--288). Our method is applicable in several contexts, in particular to flat tori and compact hyperbolic manifolds, but we will confine our discussion to the case of compact hyperbolic manifolds, since the detailed relevant background material is already available for that case in \cite{randolchapter}. The discussion of this subject in \cite{randolchapter} was motivated by analogies with classical electrostatics, but the conclusions drawn there can, if preferred, be regarded as being of a formal character, and independent of these analogies. We note that in the flat torus case the analogue of the Selberg pre-trace formula is the Poisson summation formula.

Suppose that $M$ is a compact hyperbolic manifold, and that $S$ is a finite set of distinct points in $M$. We regard the points of $S$ as mutually repelling along connecting geodesics, with the repelling magnitude varying with distance, as specified by a function $H(\rho)$ $\rho \in [0,\infty)$. In other words, the effect exerted on a point $p \in S$ by points in $S$ is computed by taking the (infinite) vector sum $\sum_n H(L_n)V_n$, where $L_n$ ranges over the lengths of the geodesic segments $g_n$ connecting $p$ to points of $S$, and $V_n$ is the unit tangent vector to $g_n$ at $p$, taken in the repelling direction. Note that $p$ is permitted to act on itself. Conceptually, this can be visualized intrinsically on $M$, or if one prefers, on the universal covering space $H^n$, where the mutual repulsions act along the countably numerous geodesics connecting points, attenuated according to $H(\rho)$. Similar considerations can also be applied to continuous, rather than point distributions on $M$, where the repulsion between infinitesimal parts of $M$ is defined analogously \cite{randolchapter}, but in this note we will be concerned with finite point distributions. In the interests of expository economy and to avoid unnecessary duplication, we refer the reader to \cite{randolchapter} for background material concerning the Selberg pretrace and trace formulas, as well as to the application of the pretrace formula to the matters considered here. In greater detail, following the discussion in \cite{randolchapter}, we assume that $H(\rho)$ is linked to a function $k(\rho)$ by the requirement that $H(\rho) = -k'(\rho)$, where $k(\rho)$, which corresponds to a point-pair invariant, is a smooth function on $[0,\infty)$ which vanishes at infinity, and for which the Selberg pretrace formula is valid, with uniform and absolute convergence. As is customary, we denote the Selberg transform of $k(\rho)$ by $h(r)$, and we will be concerned with the values of $h(r)$ on the cross-shaped subset of the complex plane composed of the union $L$ of the real axis with the closed segment of the imaginary axis from \mbox{$-(n-1)i/2$} to $(n-1)i/2$. In the standard parameterization $r \leftrightarrow \lambda$, this corresponds to eigenvalues of the positive Laplacian on $M$, which are situated on $[0,\infty)$. We will also impose, for the purposes of the present note, the additional requirement that $h(r)$ is positive on $L$, and will call the configuration $S$ stable, if the net repelling effect at each point of $S$ is zero.

Let $0=\lambda_0 < \lambda_1 \leq \lambda_2 \ldots$ be the eigenvalues of the positive Laplacian on $M$, with repetitions to account for multiplicity, and suppose $\varphi_0, \varphi_1, \varphi_2,\ldots$ is an associated orthonormal sequence of eigenfunctions. Denote the points of the set $S$ by $x_1,\ldots , x_N$. With this notation, it follows from the discussion in \cite{randolchapter} (cf.\ especially p. 287), that sets $S$ for which the quantity \[\sum_{n=1}^{\infty} \left(\sum_{i=1}^N \sum_{j=1}^N \varphi_n(x_i)\overline{\varphi_n}(x_j)\right)h(r_n)\] \begin{equation}=  \sum_{n=1}^{\infty} |\varphi_n(x_1) + \cdots \varphi_n(x_N)|^2 h(r_n) \label{condition}\end{equation} is locally minimized in the $N$-fold Cartesian product of $M$ with itself are stable in the above sense. In somewhat more detail,    formula (17) from page $283$ of \cite{randolchapter}, which describes the repelling effect produced by a point of $M$ on another point of $M$, can be used to derive this result. The requirement of zero mutual repulsion at each point of $S$ is expressible via that formula as a requirement for the simultaneous vanishing of several instances of (17). This translates, since the gradient on the $N$-fold Cartesian product of M with itself splits into the vector sum of the gradients on the factors, into a necessary condition for (1) to have a minimum.

We remark that there is a misprint on page $283$ of \cite{randolchapter}. Namely, in the formula above (17), the expression $-\nabla_x k(x,\gamma y)$ should be replaced by $-\nabla_x \sum_{\gamma} k(x,\gamma y)$, and directly above, the words ``and by the pretrace formula'' should be augmented to ``and by the pretrace formula, counting half the $r_n$'s to avoid a factor of $\frac12$.''


These observations suggest, as being of particular interest, the study of the statistical behavior, as $N \rightarrow \infty$, of point configurations that are globally minimizing. and as we will now show, such configurations are equidistributed in the limit. In order to see this, we will prove a somewhat more general fact, expressed as a theorem, from which the assertion immediately follows.

\begin{genericem}{Theorem} Suppose that $M$ is a compact Riemannian manifold,
and that $1,\varphi_1,\varphi_2,\ldots$ is a sequence of
continuous functions in $C(M)$, each having $L^2$ norm $1$, and with $\int_M \varphi_i = 0$ for $i=1,2,\ldots$. Suppose also that finite linear combinations of functions from the sequence are sup-norm dense in $C(M)$. An important special case in which these conditions are satisfied occurs when the sequence $1,\varphi_1,\varphi_2,\ldots$ is complete and orthonormal in $L^2(M)$. Denote by 
$M^N$  the $N$-fold Cartesian product of $M$ with itself. Let
$a_1,a_2,\ldots$ be a sequence of positive numbers for which the series
$\sum_{n=1}^\infty |\varphi_n(x)|^2 \, a_n$ is uniformly convergent, and
suppose, for $N=1,2,3\ldots$, that $X_N \in M^N$ globally minimizes the
quadratic expression \[ \sum_{n=1}^\infty \left| \varphi_n(x_1)+ \ldots
+\varphi_n(x_N)\right| ^2 \, a_n \,.\] Then, as $N \rightarrow \infty$,
the sets $X_N$ become equidistributed, in the sense that for any
continuous function $f$ on $M$, \[\frac{1}{N}\sum_{x \in X_N} f(x)
\rightarrow \int_M f\] as $N \rightarrow \infty$.\end{genericem}

\begin{proof} It follows from the hypotheses that for each $n$, \[
\int_{M\times\cdots\times M} \left| \varphi_n(x_1)+ \ldots +
\varphi_n(x_N)\right| ^2 \,dx_1\cdots dx_k = NV^{N-1} \,,\]where $V$ is the volume of $M$, so
\[\int_{M\times\cdots\times M} \left(\sum_{n=1}^\infty \left|
\varphi_n(x_1)+ \ldots + \varphi_n(x_N)\right| ^2 a_n\right)
\,dx_1\cdots dx_k = NV^{N-1}\sum_{n=1}^\infty a_n\,. \]

In particular, by the mean value theorem for integrals, there exists $(\overline{x}_1,\ldots ,\overline{x}_N)\in
M^N$ such that \[ \sum_{n=1}^\infty \left| \varphi_n(\overline{x}_1)+
\ldots + \varphi_n(\overline{x}_N)\right| ^2 a_n = \frac{N}{V} \sum_{n=1}^{\infty}
a_n\,,\] which shows that if $(p_1\ldots ,p_N) \in M^N$ is globally
minimizing, then \[ \sum_{n=1}^\infty \left| \varphi_n(p_1)+ \ldots +
\varphi_n(p_N)\right| ^2 a_n \leq \frac{N}{V} \sum_{n=1}^{\infty} a_n \,.\]

For a fixed $m$, we therefore conclude from the positivity of the $a_n$'s that \[ \left| \varphi_m(p_1)+
\ldots + \varphi_m(p_N)\right|^2 \leq \frac{N}{V} a_m^{-1} \sum_{n=1}^{\infty}
a_n \,,\] so \[ \left| \varphi_m(p_1)+ \ldots + \varphi_m(p_N)\right|
\leq \sqrt{\frac{N}{V} a_m^{-1} \sum_{n=1}^{\infty} a_n} \;\;,\] which implies
that \[ N^{-1} \left| \varphi_m(p_1)+ \ldots + \varphi_m(p_N)\right|
\leq \frac{C(m)}{\sqrt{N}}\,,\] where \[ C(m) = \sqrt{\frac{1}{Va_m}
\sum_{n=1}^{\infty} a_n} \;\;.\] In particular, the $X_k$'s in the limit integrate 
the uniformly dense set of functions $1,\varphi_1,\varphi_2,\ldots$ (the constant function is automatic), so Weyl's criterion is satisfied, and the $X_k$'s are equidistributed in the limit.\end{proof}

\newpage

\noindent\textbf{Remarks.}

\begin{enumerate}

\item It is not difficult to give examples of force laws for which $h$ is positive. This is exceptionally easy to see in the hyperbolic 2-dimensional case, since then $k(u) =  \frac{1}{2}\hat{h}(u)$, so there are many suitable $h,k$ pairs. In the flat torus case,  $k$ and $h$ are in all dimensions Fourier transforms of each other, which again provides many examples. We note also that although it is customary to do so, there is no formal reason to insist that the force law be always repelling.

\item If we consider globally minimizing configurations to be equivalent when they are mutual images by isometries of $M$, it would be interesting to know something about the behavior with $N$ of the number of inequivalent configurations.

\item Our definition of point interaction is intrinsic, and does not depend on an isometric embedding of $M$. A prominent example of the latter occurs in the well-known Thompson problem concerning stable configurations of electrons on the unit sphere of $R^3$, where the repelling actions resulting from the Coulomb law are transmitted via mutually connecting geodesic segments in the containing space $R^3$. Note that the Coulomb law has a singularity at the origin, whereas the methods of this paper require the function $k$ to have a smooth even extension to $R^1$, and the series which arise to be suitably convergent. This of course does not rule out the possibility of studying an intrinsic singular law as a limit of regular laws. In view of various isometric embedding theorems, there can be connections between the intrinsic and extrinsic viewpoints. For example, if distance computed in the containing space depends smoothly on intrinsic distance, the extrinsic case may in some instances become a subcase of the intrinsic case, by choosing the intrinsic repelling function $H(\rho)$ to smoothly coincide up to the diameter of $M$ with the repelling function associated with the embedding, then appropriately bringing it down to zero over a distance less than, for example, the injectivity radius of $M$, provided this can be done in such a way that the positivity requirement on $h(r)$ is satisfied. There can, however, be substantial differences between the extrinsic and intrinsic viewpoints. For example, in the case of the unit circle embedded in various ways in $R^n$, stable configurations for the embedded problem corresponding to various repelling laws in $R^n$ are sensitive to the shape of the embedded circle (cf.\ \cite{saff1}), whereas this issue is obviously not present for the intrinsic problem. Some asymptotic and extremal aspects of the intrinsic problem for geodesically equispaced points on a circle are discussed in \cite{toth1} and \cite{saff2}.

\item In the intrinsic case, it is possible for a point to act non-trivially on itself. In more detail, it follows from an examination of the analysis leading to formula (17) of \cite{randolchapter}, that zero-distance action of a point on itself plays no role, but that the effects of the non-trivial geodesics connecting a point to itself are counted. In the flat torus case, these actions cancel out in pairs, and so in that case a single point is stable with respect to itself, but in the hyperbolic case, this is not generally true. In the hyperbolic case, the problem becomes one of studying minima for the function \[ \sum_{n=1}^{\infty} |\varphi_n(x)|^2 h(r_n)\] on $M$.

\end{enumerate}

\newcommand{\noopsort}[1]{}

\bigskip

\begin{flushleft}
{\sc Ph.D Program in Mathematics\\CUNY Graduate Center\\365 Fifth
Avenue\\New York, NY 10016}
\end{flushleft}

\end{document}